\documentclass[12pt,a4paper]{article}
\usepackage{amsmath,amssymb,amsbsy,amscd,amsthm}
\newtheorem{thm}{Theorem}

\newcommand{\nd }{\noindent}
\newcommand{\TA }{{\rm TA}}
\newcommand{\GA }{{\rm GA}}
\newcommand{\ch }{\mathop{\rm char}\nolimits}
\newcommand{\zs}{\{ 0\} }
\newcommand{\sm}{\setminus}
\newcommand{\C}{{\bf C}}
\newcommand{\Q}{{\bf Q}}

\newcommand{\x}{{\bf x}}
\newcommand{\Rx}{R[{\bf x}]}
\newcommand{\ep}{{\epsilon}}
\newcommand{\p}{{\mathfrak{p}}}
\newcommand{\Spec }{\mathop{\rm Spec}\nolimits}

\begin{document}
\title{Degeneration of tame automorphisms of a polynomial ring}

\author{Shigeru Kuroda
\thanks{Partly supported by the Grant-in-Aid for 
Young Scientists (B) 24740022, 
Japan Society for the Promotion of Science. }}


\date{}

\maketitle

\begin{abstract}
Recently, 
Edo-Poloni constructed a family of tame automorphisms 
of a polynomial ring in three variables 
which degenerates to a wild one. 
In this note, 
we generalize the example by a different method. 
\end{abstract}

\section{Introduction}\label{sect:intro}
\setcounter{equation}{0}

For each commutative ring $R$, 
let $\Rx =R[x_1,\ldots ,x_n]$ 
be the polynomial ring in $n$ variables over $R$, 
and $\GA _n(R)$ 
the automorphism group of the $R$-algebra $\Rx $. 
We express $\phi \in \GA _n(R)$ as 
$(\phi (x_1),\ldots ,\phi (x_n))\in \Rx ^n$. 
We say that $\phi \in \GA _n(R)$ is {\it affine} 
if the total degrees of 
$\phi (x_1),\ldots ,\phi (x_n)$ are equal to one, 
{\it triangular} 
if $\phi (x_i)$ belongs to 
$R^*x_i+R[x_1,\ldots ,x_{i-1}]$ 
for each $i$, 
and {\it elementary} if 
there exist $1\leq l\leq n$ 
and $q\in R[\{ x_i\mid i\neq l\} ]$ 
such that $\phi (x_l)=x_l+q$ and 
$\phi (x_i)=x_i$ for $i\neq l$. 
Let $\TA _n(R)$ be the subgroup of $\GA _n(R)$ 
generated by affine automorphisms 
and triangular automorphisms. 
We say that $\phi \in \GA _n(R)$ is {\it tame} 
if $\phi $ belongs to $\TA _n(R)$, 
and {\it wild} otherwise. 
If $K$ is a field, 
then $\TA _2(K)$ is equal to $\GA _2(K)$ 
by Jung~\cite{Jung} and van der Kulk~\cite{Kulk}. 
In 2004, 
Shestakov-Umirbaev~\cite{SU} 
showed that the same does not hold when 
$\ch K=0$ and $n=3$. 
Today, 
a number of elements of $\GA _3(K)$ 
are known to be wild 
thanks to the Shestakov-Umirbaev theory 
and its modification~\cite{tame3}.

Recently, 
Edo-Poloni~\cite[\S 4]{EP} 
constructed 
$\phi \in \GA _3(\C [t])$ 
with the following properties, 
where $t$ is a variable, 
and $(x_1,x_2,x_3)=(z,y,x)$ 
in their notation: 

\smallskip 

\nd (A) 
$\phi $ factors as 
$\tau ^{-1}\circ \ep \circ \tau $ 
in $\GA _3(\C [t,t^{-1}])$, 
where $\tau \in \GA _3(\C [t,t^{-1}])$ 
is triangular 
and $\ep \in \GA _3(\C [t])$ is elementary. 

\nd (B) 
$\phi _0$ is the wild automorphism 
$\exp x_1^{2l}(x_1x_3+x_2^{l+1})\delta _0$, 
where $l\geq 1$ 
and $\delta _0$ is the 
$\C [x_1]$-derivation of $\C [\x ]$ defined by 
$\delta _0(x_2)=x_1$ and 
$\delta _0(x_3)=-(l+1)x_2^l$.

\smallskip

Here, 
for each $\alpha \in \C $, 
let $\phi _{\alpha }$ be the automorphism of 
$\C [\x ]=(\C [t]/(t-\alpha ))\otimes _{\C [t]}\C [t][\x ]$ 
induced by $\phi $. 
From (A), 
it follows that $\phi _{\alpha }$ 
is tame if $\alpha \neq 0$. 
Thus, 
a wild automorphism is obtained as a ``limit" 
of tame automorphisms. 
Using this example, 
Edo-Poloni~\cite[Cor.~4.3]{EP} concluded that 
$\TA _3(\C )$ is not closed in 
the ind-group $\GA _3(\C )$. 
The purpose of this note 
is to generalize Edo-Poloni's example 
by a different method. 
From our construction, 
we easily see why such a phenomenon occurs.

The first version of this note was written for 
a discussion 
at Saitama University in 2014. 
Thanks are due to the participants, 
especially to 
Prof.~E. Edo for introducing this topic, 
and to Prof.~A. Dubouloz for the suggestion 
of generalizing $R$ 
which was originally 
a multivariate polynomial ring.

\section{Construction}\label{sect:construct}
\setcounter{equation}{0}

Let $R$ be any commutative $\Q $-algebra, 
and $\delta $ a {\it triangular} $R$-derivation 
of $\Rx $, 
i.e., 
an $R$-derivation of $\Rx $ such that 
$f_i:=\delta (x_i)$ belongs to $R[x_1,\ldots ,x_{i-1}]$ 
for each $i$. 
Then, 
for each $h\in \Rx $ with $\delta (h)=0$, 
we can define 
$\phi :=\exp h\delta \in \GA _n(R)$ 
by $(\exp h\delta )(q)=\sum _{l\geq 0}h^l\delta ^l(q)/l!$ 
for $q\in \Rx $. 
For each $\p \in \Spec R$, 
let $\phi _{\p }$ be 
the element of $\GA_n(\kappa (\p ))$ 
induced by $\phi $, 
where $\kappa (\p )$ is the residue field 
of the localization $R_{\p }$. 
We consider when $\phi _{\p }$ is tame.

First, 
recall that, 
when $n=3$, $R$ is a field and $f_1=0$, 
we have $\phi \not\in \TA _3(R)$ 
if and only if 
$f_2\neq 0$, 
$h\not\in R[x_1]$ and 
$\partial _{x_2}(f_3)\not\in f_2R[x_1,x_2]$ 
(cf.~\cite[Thm.~3.2.3 (and the preceding remark)]{wild3}). 
Hence, 
we get the following theorem, 
where 
$\bar{f}$ denotes the image of $f$ in 
$\kappa (\p )[\x ]$ for each $f\in \Rx $.

\begin{thm}\label{thm:main}
Assume that $n=3$ and $f_1\in \p $. 
Then, 
$\phi _{\p }$ is wild 
if and only if 
$\bar{f}_2\neq 0$, 
$\bar{h}\not\in \kappa (\p )[x_1]$ 
and 
$\partial _{x_2}(\bar{f}_3)\not\in 
\bar{f}_2\kappa (\p )[x_1,x_2]$. 
\end{thm}

To discuss the case where $f_1\not\in \p $, 
assume that 
$f_1$ is not a nilpotent element of $R$, 
and let $R'$ be the localization $R_{f_1}$. 
Since no confusion arises, 
we use the letter $\phi $ to denote 
the element of $\GA _n(R')$ induced by $\phi $, 
and $\delta $ to denote 
the derivation of $R'[\x ]$ 
induced by $\delta $, 
whose kernel is denoted by $R'[\x ]^{\delta }$. 
Then, 
the image of $h$ in $R'[\x ]$ 
belongs to $R'[\x ]^{\delta }$. 
In this situation, 
there exists a triangular automorphism 
$\tau =(x_1,g_2,\ldots ,g_n)$ of $R'[\x ]$ 
such that $R'[\x ]^{\delta }=R'[g_2,\ldots ,g_n]$. 
Actually, 
since $\delta (x_1/f_1)=1$, 
we have $R'[\x ]=R'[\x ]^{\delta }[x_1/f_1]$, 
and 
$$
\sigma :R'[\x ]\ni 
q\mapsto \sum _{l\geq 0}
\frac{\delta ^l(q)}{l!}
(-x_1/f_1)^l
\in R'[\x ]
$$
is a homomorphism 
of $R'$-algebras 
satisfying $\sigma (R'[\x ])=R'[\x ]^{\delta }$ 
and $\sigma (x_1)=0$ 
(cf.~e.g.~\cite[1.3.21 and 1.3.23]{Essen}). 
Now, 
let $p\in R'[x_2,\ldots ,x_n]$ 
be such that $\tau (p)=h$ in $R'[\x ]$. 
Then, in $\GA _n(R')$, 
we have 
\begin{equation}\label{eq:A}
\begin{aligned}
&\tau ^{-1}\circ \phi \circ \tau 
=\exp 
\bigl(\tau ^{-1}\circ (h\delta )\circ \tau \bigr) 
=\exp 
\bigl(\tau ^{-1}\circ (\tau (p)\delta )
\circ \tau \bigr) \\ 
&\quad =\exp (p\tau ^{-1}\circ \delta 
\circ \tau ) =\exp pf_1\partial _{x_1}
=(x_1+f_1p,x_2,\ldots ,x_n)=:\ep , 
\end{aligned}\tag{$*$}
\end{equation}
and so $\phi =\tau \circ \ep \circ \tau ^{-1}$. 
Therefore, 
the following theorem 
holds for any $n$.

\begin{thm}\label{thm:tame}
If $f_1\not\in \p $, then $\phi _{\p }$ is tame. 
\end{thm}

It is interesting to note that 
the extension 
$\tilde{\phi }\in \GA _{n+1}(R)$ 
of $\phi $ defined by $\tilde{\phi }(x_{n+1})=x_{n+1}$ 
is tame by Smith~\cite{Smith}. 
In fact, 
let $\tilde{\delta }$ be the extension of $\delta $ 
to $R[x_1,\ldots ,x_{n+1}]$ 
defined by $\tilde{\delta }(x_{n+1})=0$, 
and 
$\gamma $ the elementary automorphism 
$(x_1,\ldots ,x_n,x_{n+1}+h)$. 
Then, 
$\rho :=\exp x_{n+1}\tilde{\delta }$ is tame and 
$\tilde{\phi }
=\gamma ^{-1}\circ \rho ^{-1}\circ \gamma \circ \rho $.

Finally, 
we construct 
$\phi \in \GA _3(\C [t])$ 
satisfying (A) and (B) 
using our method. 
Let $R=\C [t]$, 
$f_1=t$, $f_2=x_1$ and $f_3=-(l+1)x_2^l$, 
where $n=3$. 
Then, 
we have $R'=\C [t,t^{-1}]$. 
Observe that $\delta $ kills 
$$
g_2:=x_2-x_1^2/(2t)
\quad\text{and}\quad
g_3:=x_3+\sum _{i=0}^lc_it^{-(i+1)}x_1^{2i+1}x_2^{l-i}, 
$$
and that $\tau :=(x_1,g_2,g_3)\in \GA _3(\C [t,t^{-1}])$ 
is triangular, 
where $c_0=l+1$ and 
$c_1\ldots ,c_l\in \Q $ are defined by 
$c_i(2i+1)=-c_{i-1}(l-i+1)$ by induction. 
Set 
$$
p=(c_lt^l/2)\left((2x_2)^{2l+1}+t(x_3/c_l)^2\right) 
\quad\text{and}\quad
\ep =(x_1+tp,x_2,x_3). 
$$
Note that 
$h:=\tau (p)$ is killed by $\delta $. 
In the following, 
we check $h\in \C[t][\x ]$ 
and $h|_{t=0}=x_1^{2l}(x_1x_3+x_2^{l+1})$. 
Here, 
for each $q\in \C [t][\x ]$, 
we denote by $q|_{t=0}$ the element of $\C [\x ]$ 
obtained from $q$ by the substitution $t\mapsto 0$. 
Then, 
it follows that 
$\phi =\exp h\delta $ satisfies (B). 
Since 
$\phi =\tau \circ \ep \circ \tau ^{-1}$ 
by (\ref{eq:A}), 
(A) is also satisfied.

We may write 
$h=t^{l+1}x_3^2/(2c_l)+x_1^{2l+1}x_3+q$, 
where 
$q\in x_2\C [t,t^{-1}][x_1,x_2]+tx_3\C [t][x_1,x_2]$. 
Since the monomial $x_1^{2l+1}x_3$ appears in $h$, 
the minimal integer $r$ 
for which 
$t^rh$ belongs to $\C [t][\x ]$ 
is nonnegative. 
If $r$ is positive, 
then $h':=(t^rh)|_{t=0}=(t^rq)|_{t=0}$ 
belongs to $x_2\C [x_1,x_2]\sm \zs $. 
Let $\delta _0$ be 
the derivation of $\C [\x ]$ as in (B). 
Then, 
we have $\delta _0(h')=x_1\partial _{x_2}(h')\neq 0$, 
while $\delta _0(h')=\delta (t^rh)|_{t=0}=0$ 
since $\delta (t^rh)=0$. 
This is a contradiction. 
Thus, 
we get $r=0$, 
and so $h'=h|_{t=0}=x_1^{2l+1}x_3+q|_{t=0}$. 
Since $\delta _0(h')=0$ 
and 
$q|_{t=0}\in x_2\C [x_1,x_2]$, 
it follows that 
$q|_{t=0}=x_1^{2l}x_2^{l+1}$, 
proving the claim.

\bigskip 

\nd 
Note: 
Recall that the {\it length} $\lambda (\theta )$ 
of a tame automorphism $\theta $ 
is by definition the minimal number of 
triangular automorphisms 
needed to express $\theta $ 
together with affine automorphisms. 
Due to Furter~\cite{Furter}, 
$\lambda $ is a lower semi-continuous 
function on $\GA _2(\C)=\TA _2(\C )$. 
So, 
if $\psi \in \GA _2(\C [t])$ 
and $l\geq 0$ are such that 
$\lambda (\psi _{\alpha })\leq l$ 
for all $\alpha \in \C ^*$, 
then we have $\lambda (\psi _0)\leq l$. 
Edo-Poloni~\cite[\S 1]{EP} remarked that 
this fails if $n=3$ because of their example. 
Edo raised a question whether 
there exists a counterexample 
in which $\psi _0$ is also tame. 
We claim that the answer is yes if $n=4$: 
Let $\psi \in \GA_4(\C [t])$ 
be the extension of the above $\phi $ 
defined by $\psi (x_4)=x_4$. 
Then, 
we have 
$\lambda (\psi _{\alpha })\leq 3$ 
for all $\alpha \in \C ^*$, 
while $\lambda (\psi _0)\leq 4$ 
by Smith~\cite{Smith}. 
Moreover, 
we can show that 
$\lambda (\psi _0)$ is in fact 
equal to four 
(the details will appear in a future paper).

\noindent
Department of Mathematics and Information Sciences\\ 
Tokyo Metropolitan University \\
1-1  Minami-Osawa, Hachioji \\
Tokyo 192-0397, Japan\\
kuroda@tmu.ac.jp


\begin{thebibliography}{00}


\bibitem{EP}
E. Edo and P.-M. Poloni, 
On the closure of the tame automorphism group 
of affine three-space, 
arXiv:math.AG/1403.2843v2. 



\bibitem{Essen}A.~van den Essen, 
{\it Polynomial automorphisms and the Jacobian conjecture}, 
Progress in Mathematics, 
Vol.\ 190, Birkh\"auser, Basel, Boston, Berlin, 
2000. 


\bibitem{Furter}
J.-P. Furter, 
On the length of polynomial automorphisms of the affine plane, 
Math. Ann. {\bf 322} (2002), 
401--411. 


\bibitem{Jung}H.~Jung, 
\"Uber ganze birationale Transformationen der Ebene, 
J.\ Reine Angew.\ Math.\ {\bf 184} (1942), 161--174. 


\bibitem{Kulk}W.~van der Kulk, 
On polynomial rings in two variables, 
Nieuw Arch.\ Wisk. (3) {\bf 1} (1953), 33--41. 


\bibitem{tame3}S. Kuroda, 
Shestakov-Umirbaev reductions 
and Nagata's conjecture on a polynomial automorphism, 
Tohoku Math.\ J.\ {\bf 62} (2010), 75--115.


\bibitem{wild3}S.~Kuroda, 
Wildness of polynomial automorphisms in three variables, 
arXiv:math.AC/1110.1466v1. 


\bibitem{SU}
I.~Shestakov and U.~Umirbaev, 
The tame and the wild automorphisms of polynomial rings 
in three variables, 
J.\ Amer.\ Math.\ Soc.\ {\bf 17} (2004), 197--227. 

\bibitem{Smith}M. K. Smith, 
Stably tame automorphisms, 
J. Pure Appl. Algebra {\bf 58} (1989), 
209--212. 

\end{thebibliography}
\end{document}